\input amstex

\documentstyle{amsppt}
  \magnification=1100
  \hsize=6.2truein
  \vsize=9.0truein
  \hoffset 0.1truein
  \parindent=2em

\NoBlackBoxes


\font\eusm=eusm10                   


\font\eusms=eusm7                       

\font\eusmss=eusm5                      


\newcount\theTime
\newcount\theHour
\newcount\theMinute
\newcount\theMinuteTens
\newcount\theScratch
\theTime=\number\time
\theHour=\theTime
\divide\theHour by 60
\theScratch=\theHour
\multiply\theScratch by 60
\theMinute=\theTime
\advance\theMinute by -\theScratch
\theMinuteTens=\theMinute
\divide\theMinuteTens by 10
\theScratch=\theMinuteTens
\multiply\theScratch by 10
\advance\theMinute by -\theScratch
\def\timeHHMM{{\number\theHour:\number\theMinuteTens\number\theMinute}}

\def\today{{\number\day\space
 \ifcase\month\or
  January\or February\or March\or April\or May\or June\or
  July\or August\or September\or October\or November\or December\fi
 \space\number\year}}

\def\timeanddate{{\timeHHMM\space o'clock, \today}}

\define\Afr{{\frak A}}

\define\Afrinf{{\Afr_{(-\infty,\infty)}}}

\define\Afrt{{\widetilde\Afr}}

\define\Afrtinf{{\Afrt_{(-\infty,\infty)}}}

\define\ah{{\hat a}}

\define\Ao{{A\oup}}

\define\Bfr{{\frak B}}

\define\biggnm#1{
  \bigg|\bigg|#1\bigg|\bigg|}

\define\bignm#1{
  \big|\big|#1\big|\big|}

\define\bh{\hat b}

\define\Bo{{B\oup}}

\define\clspan{\overline\lspan}

\define\Cpx{\bold C}

\define\dif{\text{\it d}}

\define\eqdef{{\;\overset\text{def}\to=\;}}

\define\Eto#1{E_{(\to{#1})}}

\define\fpamalg#1{{\dsize\;\operatornamewithlimits*_{#1}\;}}

\define\fpiamalg#1{{\tsize\;({*_{#1}})_{\raise-.5ex\hbox{$\ssize\iota\in I$}}}}

\define\freeprod#1#2{\mathchoice
     {\operatornamewithlimits{\ast}_{#1}^{#2}}
     {\raise.5ex\hbox{$\dsize\operatornamewithlimits{\ast}
      _{#1}^{#2}$}\,}
     {\text{oops!}}{\text{oops!}}}

\define\freeprodi{\mathchoice
     {\operatornamewithlimits{\ast}
      _{\iota\in I}}
     {\raise.5ex\hbox{$\dsize\operatornamewithlimits{\ast}
      _{\sssize\iota\in I}$}\,}
     {\text{oops!}}{\text{oops!}}}

\define\freeprodvni{\mathchoice
      {\operatornamewithlimits{\overline{\ast}}
       _{\iota\in I}}
      {\raise.5ex\hbox{$\dsize\operatornamewithlimits{\overline{\ast}}
       _{\sssize\iota\in I}$}\,}
      {\text{oops!}}{\text{oops!}}}

\define\GNS{{\text{\rm GNS}}}

\define\Hil{{\mathchoice
     {\text{\eusm H}}
     {\text{\eusm H}}
     {\text{\eusms H}}
     {\text{\eusmss H}}}}


\define\Hilo{\Hil\oup}

\define\Hilto#1{\Hil_{(\to{#1})}}

\define\id{\text{\rm id}}

\define\Integers{\bold Z}

\define\Lambdao{{\Lambda\oup}}

\define\ld#1{{\hbox{..}(#1)\hbox{..}}}

\define\lrnm#1{\left\|#1\right\|}

\define\lspan{\text{\rm span}@,@,@,}

\define\nm#1{\|#1\|}

\define\Naturals{{\bold N}}

\define\otdts#1{\otimes_{#1}\cdots\otimes_{#1}}

\define\oup{^{\text{\rm o}}}

\define\owedge{{
     \operatorname{\raise.5ex\hbox{\text{$
     \ssize{\,\bigcirc\llap{$\ssize\wedge\,$}\,}$}}}}}

\define\owedgeo#1{{
     \underset{\raise.5ex\hbox
     {\text{$\ssize#1$}}}\to\owedge}}

\define\Pto#1{{P_{(\to{#1})}}}


\define\pup#1#2{{{\vphantom{#2}}^{#1}\!{#2}}\vphantom{#2}}

\define\QED{\newline
            \line{$\hfill$\qed}\enddemo}

\define\restrict{\lower .3ex
     \hbox{\text{$|$}}}

\define\smd#1#2{\underset{#2}\to{#1}}

\define\smdb#1#2{\undersetbrace{#2}\to{#1}}

\define\smdbp#1#2#3{\overset{#3}\to
     {\smd{#1}{#2}}}

\define\smdbpb#1#2#3{\oversetbrace{#3}\to
     {\smdb{#1}{#2}}}

\define\smdp#1#2#3{\overset{#3}\to
     {\smd{#1}{#2}}}

\define\smdpb#1#2#3{\oversetbrace{#3}\to
     {\smd{#1}{#2}}}

\define\smp#1#2{\overset{#2}\to
     {#1}}

\define\Thetainf{{\Theta_{(-\infty,\infty)}}}

\define\tint{{\tsize\int}}

\define\tocdots
  {\leaders\hbox to 1em{\hss.\hss}\hfill}

\define\Tr{\text{\rm Tr}}


  \newcount\mycitestyle \mycitestyle=1 

  \newcount\bibno \bibno=0
  \def\newbib#1{\advance\bibno by 1 \edef#1{\number\bibno}}
  \ifnum\mycitestyle=1 \def\cite#1{{\rm[\bf #1\rm]}} \fi
  \def\scite#1#2{{\rm[\bf #1\rm, #2]}}


  \newcount\ignorsec \ignorsec=0
  \def\notasec{\ignorsec=1}

  \newcount\secno \secno=0
  \def\newsec#1{\procno=0 \subsecno=0 \ignorsec=0
    \advance\secno by 1 \edef#1{\number\secno}
    \edef\currentsec{\number\secno}}

  \newcount\subsecno
  \def\newsubsec#1{\procno=0 \advance\subsecno by 1
    \edef\currentsec{\number\secno.\number\subsecno}
     \edef#1{\currentsec}}

  \newcount\appendixno \appendixno=0
  \def\newappendix#1{\procno=0 \ignorsec=0 \advance\appendixno by 1
    \ifnum\appendixno=1 \edef\appendixalpha{\hbox{A}}
      \else \ifnum\appendixno=2 \edef\appendixalpha{\hbox{B}} \fi
      \else \ifnum\appendixno=3 \edef\appendixalpha{\hbox{C}} \fi
      \else \ifnum\appendixno=4 \edef\appendixalpha{\hbox{D}} \fi
      \else \ifnum\appendixno=5 \edef\appendixalpha{\hbox{E}} \fi
      \else \ifnum\appendixno=6 \edef\appendixalpha{\hbox{F}} \fi
    \fi
    \edef#1{\appendixalpha}
    \edef\currentsec{\appendixalpha}}

  \newcount\procno \procno=0
  \def\newproc#1{\advance\procno by 1
   \ifnum\ignorsec=0 \edef#1{\currentsec.\number\procno}
                     \edef\currentproc{\currentsec.\number\procno}
   \else \edef#1{\number\procno}
         \edef\currentproc{\number\procno}
   \fi}

  \newcount\subprocno \subprocno=0
  \def\newsubproc#1{\advance\subprocno by 1
   \ifnum\subprocno=1 \edef#1{\currentproc a} \fi
   \ifnum\subprocno=2 \edef#1{\currentproc b} \fi
   \ifnum\subprocno=3 \edef#1{\currentproc c} \fi
   \ifnum\subprocno=4 \edef#1{\currentproc d} \fi
   \ifnum\subprocno=5 \edef#1{\currentproc e} \fi
   \ifnum\subprocno=6 \edef#1{\currentproc f} \fi
   \ifnum\subprocno=7 \edef#1{\currentproc g} \fi
   \ifnum\subprocno=8 \edef#1{\currentproc h} \fi
   \ifnum\subprocno=9 \edef#1{\currentproc i} \fi
   \ifnum\subprocno>9 \edef#1{TOO MANY SUBPROCS} \fi
  }

  \newcount\tagno \tagno=0
  \def\newtag#1{\advance\tagno by 1 \edef#1{\number\tagno}}



\notasec
  \newtag{\cstafp}
\newsec{\Notation}
 \newproc{\notationAo}
 \newproc{\notationOminus}
 \newproc{\notationLambdao}
\newsec{\CondExp}
 \newproc{\ProjGivesExp}
 \newproc{\Eominus}
  \newtag{\EAominusB}
  \newtag{\AAominusBB}
 \newproc{\FreeProdOfCondExp}
 \newproc{\FPCondExp}
 \newproc{\ProjOntoFactor}
\newsec{\PIFP}
 \newproc{\mostgeneral}
   \newtag{\phiofpn}
   \newtag{\xinTheta}
  \newsubproc{\alphaouter}
  \newsubproc{\phiThetao}
   \newtag{\LamLamkerphi}
  \newsubproc{\AfrFreenegone}
   \newtag{\heresx}
   \newtag{\inThetaneg}
  \newsubproc{\AfrZeroSimple}
  \newsubproc{\AfrNegSimple}
  \newsubproc{\pnkFree}
   \newtag{\Eofpm}
  \newsubproc{\HeredProj}
 \newproc{\examplesPI}

\newbib{\Avitzour}
\newbib{\CuntzZZAddMult}
\newbib{\CuntzZZDimFunct}
\newbib{\DykemaZZFPFDNT}
\newbib{\DykemaZZFaithful}
\newbib{\DykemaZZSimplicity}
\newbib{\DykemaHaagerupRordam}
\newbib{\DykemaRordamZZPI}
\newbib{\DykemaRordamZZProj}
\newbib{\RordamZZUHFII}
\newbib{\Tomiyama}
\newbib{\VoiculescuZZSymmetries}
\newbib{\VDNbook}

\topmatter
  \title Purely Infinite, Simple $C^*$-algebras Arising from Free
         Product Constructions, II
  \endtitle

  \author Kenneth J\. Dykema \endauthor

  \date 22 August, 1999 \enddate

  \rightheadtext{Purely infinite, simple $C^*$--algebras, II, \timeanddate}

  \leftheadtext{Purely infinite, simple $C^*$--algebras, II, \timeanddate}

  \address Department of Mathematics, Texas A\&M University, College Station TX 77843, USA \endaddress
  \email {\tt Ken.Dykema\@math.tamu.edu,}\hskip.3em
         {\it Internet:} {\tt http://www.math.tamu.edu/\~{\hskip0.1em}Ken.Dykema/}
  \endemail

  \abstract
     Certain reduced free products of C$^*$--algebras with respect to faithful
     states are simple and purely infinite.
  \endabstract

  \subjclass 46L05, 46L35 \endsubjclass

\endtopmatter

\document \TagsOnRight \baselineskip=18pt

\heading Introduction.\endheading
\vskip3ex

Given unital C$^*$--algebras $A$ and $B$ with states $\phi_A$ and $\phi_B$,
whose GNS representation are faithful,
their reduced free product C$^*$--algebra,
$$ (\Afr,\phi)=(A,\phi_A)*(B,\phi_B), \tag{\cstafp} $$
was introduced in~\cite{\VoiculescuZZSymmetries} and~\cite{\Avitzour}.
It is the natural construction in Voiculescu's free probability theory
(see~\cite{\VDNbook}), and
Voiculescu's theory has been vital to the study of
these C$^*$--algebras.

In~\cite{\Avitzour}, D\. Avitzour showed that reduced free product
C$^*$--algebra $\Afr$ in~(\cstafp) is simple if $A$ and $B$ are not too small
(with respect to their states), in that they have enough orthogonal unitaries.

A unital C$^*$--algebra is said to be {\it infinite} if it contains a
nonunitary isometry, and is said to be {\it purely infinite} if every
hereditary C$^*$--subalgebra of it is infinite.
One of the most interesting open questions about simple C$^*$--algebras
is whether every infinite simple C$^*$--algebra is purely infinite.
In~\cite{\DykemaRordamZZProj}, M\.~R\o{}rdam and the author showed that in the
free product~(\cstafp), if  $\phi_A$ and $\phi_B$ are faithful, if one of them
is not a trace and if $A$ and $B$ are not too small in that they satisfy
a condition like (but slightly weaker than) Avitzour's condition, then the free
product
C$^*$--algebra $\Afr$ is infinite, and furthermore, it is properly infinite.
It remained open whether these C$^*$--algebras are purely infinite, or indeed
whether any C$^*$--algebras arising as reduced free products with respect to
faithful states are purely infinite.
(In~\cite{\DykemaRordamZZPI}, the same authors had shown that some
C$^*$--algebras arising as reduced free
product with respect to nonfaithful states are purely infinite.)

In this paper, we show that certain free product C$^*$--algebras $\Afr$
in~(\cstafp), with respect to faithful states, are simple and purely infinite.
For example, if $A=B=M_2(\Cpx)$ and if $\phi_A$ and $\phi_B$ are faithful states
on $M_2(\Cpx)$ that are not unitarily equivalent then $\Afr$ is simple and
purely infinite.
Note that it was previously not known whether any of these particular reduced
free product C$^*$--algebras of $M_2(\Cpx)$ with $M_2(\Cpx)$ were even simple.
(See Examples~\examplesPI{} for some more examples.)

The most striking condition that we require of $(A,\phi_A)$ and $(B,\phi_B)$
for the free product
C$^*$--algebra~(\cstafp) to be purely infinite and simple is that one of the
algebras, say $A$, contain a partial isometry, $v$, with orthogonal domain and
range projections and scaling the state $\phi_A$ by $\lambda$ for some
$0<\lambda<1$, i.e\. that $\phi_A(va)=\lambda\phi_A(av)$ for every $a\in A$.
These are fairly strong conditions, but they do arise in a large number of
situations.
Although both the statement and proof of the main result, Theorem~\mostgeneral,
are quite technical, we believe they comprise an important advance in the
understanding of simplicity and infiniteness of reduced free product
C$^*$--algebras.

\vskip3ex
\heading\S\Notation.  Notation.\endheading
\vskip3ex

Let us briefly describe some notation.

\proclaim{\notationAo} \rm
Given a unital C$^*$--algebra $A$ and a state, $\psi$ of $A$, (which will
usually be implicit from the context), we use the symbol $\Ao$ to denote the
kernel of $\psi$.
\endproclaim

\proclaim{\notationOminus} \rm
If $A$ is a C$^*$--algebra with a state $\phi$, then for any
C$^*$--subalgebra, $D\subseteq A$, we define
$$ A\ominus D=\{a\in A\mid\forall d\in D,\,\phi(da)=0\}. $$
This is just the orthocomplement of $D$ in the Hilbert space of the GNS
representation, pulled back to $A$.
\endproclaim

\proclaim{\notationLambdao}\rm
If $X_\iota$ ($\iota\in I$) are subsets of an algebra $A$, then
$$ \Lambdao\bigl((X_\iota)_{\iota\in I}\bigr)
=\{x_1x_2\cdots x_n\mid n\in\Naturals,\,x_j\in X_{\iota_j},\,
\iota_1\neq\iota_2,\iota_2\neq\iota_3,\ldots,\iota_{n-1}\neq\iota_n\}, $$
(written simply $\Lambdao(X_1,X_2)$ if $I=\{1,2\}$), is the set of all
``alternating words'' in the $X_\iota$.
\endproclaim

\vskip3ex
\heading\S\CondExp.  Conditional expectations.\endheading
\vskip3ex

Tomiyama~\cite{\Tomiyama} proved that if $B\subseteq A$ is a C$^*$--subalgebra
and if $E:A\to B$ is a projection of norm~$1$ then $E$ is positive and a
conditional expectation.
In this section, we consider the reduced free product of conditional
expectations, a special case of which will be important in the sequel.
The results proved here are perhaps well known, but their proofs
are included out of a desire for completeness.
This section can be viewed as a C$^*$--version
of~\scite{\DykemaZZFPFDNT}{\S3}.

\proclaim{Lemma \ProjGivesExp}
Let $A$ be a unital C$^*$--algebra and let $B\subseteq A$ be a unital
C$^*$--subalgebra.
Let $\phi$ be a faithful state on $A$ and let $(\pi,\Hil,\xi)=\GNS(A,\phi)$, so
that $\{\ah\mid a\in A\}$ is a dense subspace of $\Hil$ with inner product
$\langle\ah_1,\ah_2\rangle=\phi(a_2^*a_1)$.
Let $\Hil_B=\overline{\{\bh\mid b\in B\}}\subseteq\Hil$ and let
$P$ be the projection from $\Hil$ onto $\Hil_B$.
Suppose that for some norm--dense subset, $X$, of $A$, we have
$P\pi(x)\restrict_{\Hil_B}\in P\pi(B)\restrict_{\Hil_B}$ for every $x\in X$.
Then there is a projection, $E$, of norm~$1$ from $A$ onto $B$, satisfying
$\phi\circ E=\phi$.
\endproclaim
\demo{Proof}
Since $\phi$ is faithful, $\pi$ is a faithful representation of $A$.
Let $\pi_B:B\to\text{\eusm B}(\Hil_B)$ be
$\pi_B(b)=P\pi(b)\restrict_{\Hil_B}$.
Then $(\pi_B,\Hil_B,\xi)=\GNS(B,\phi\restrict_B)$, and since $\phi\restrict_B$ is faithful on $B$,
$\pi_B$ is a faithful representation of $B$.
Now from $P\pi(X)\restrict_{\Hil_B}\subseteq\pi_B(B)$ and taking limits in norm
we get that $P\pi(A)\restrict_{\Hil_B}\subseteq\pi_B(B)$, and hence we may
define $E:A\to B$ by $E(a)=\pi_B^{-1}(P\pi(a)\restrict_{\Hil_B})$.
Then clearly $E(b)=b$ if $b\in B$ and $\nm{E(a)}\le\nm a$ for every $a\in A$,
so $E$ is a projection of norm~$1$.
Moreover,
$$ \phi(E(a))=\langle \pi_B(E(a))\xi,\xi\rangle
=\langle P\pi(a)\xi,\xi\rangle
=\langle\pi(a)\xi,\xi\rangle=\phi(a), $$
so $E$ preserves $\phi$.
\QED

Let us here recall (see Lemma~3.2 of~\cite{\DykemaZZFPFDNT}) that the converse
of the above lemma is true, namely, if $B$ is a C$^*$--subalgebra of the
C$^*$--algebra $A$, and if $E:A\to B$ is a projection of norm~$1$ such that
$\phi\circ E=\phi$ for a faithful state $\phi$, on $A$, then $E$ is implemented
by a projection $P$ in Hilbert space of the GNS representation.

\proclaim{Lemma \Eominus}
Let $A$ be a C$^*$--algebra with a faithful state $\phi$ and let $E:A\to B$ be
a projection of norm~$1$ onto a C$^*$--subalgebra, $B$, of $A$.
Suppose $\phi\circ E=\phi$.
Then
$$ \ker E=A\ominus B \tag{\EAominusB} $$
and
$$ A=(A\ominus B)+B. \tag{\AAominusBB} $$
\endproclaim
\demo{Proof}
Let $(\pi,\Hil,\xi)=\GNS(A,\phi)$ and let $P$ be the projection from $\Hil$
onto $\Hil_B=\{\bh\mid b\in B\}$.
Then by~\scite{\DykemaZZFPFDNT}{3.2},
$E(a)\hat{\;}=P\ah$.
This shows that $\ker E\subseteq A\ominus B$.
To show the opposite inclusion, suppose $a\in A\ominus B$.
Then $0=E(E(a)^*a)=E(a)^*E(a)$, so $E(a)=0$.
Now~(\AAominusBB) follows from~(\EAominusB).
\QED

\proclaim{Proposition \FreeProdOfCondExp}
Suppose $A_\iota$ are unital C$^*$--algebras with faithful
states $\phi_\iota$, for $\iota$ in some index set $I$.
Suppose $B_\iota\subseteq A_\iota$ are unital C$^*$--subalgebras with
projections of norm $1$, $E_\iota:A_\iota\rightarrow B_\iota$, such that
$\phi_\iota\circ E_\iota=\phi_\iota$.
Let $(\Afr,\phi)=\freeprodi(A_\iota,\phi_\iota)$ be the free product of
C$^*$--algebras and consider the
C$^*$--subalgebra $\Bfr=C^*(\bigcup_{\iota\in I}B_\iota)\subseteq\Afr$.
Then there is a projection of norm $1$, $E:\Afr\rightarrow\Bfr$ such that
$\phi\circ E=\phi$
and $E(a)=E_\iota(a)$ whenever $a\in A_\iota$.
\endproclaim
\demo{Proof}
Let $(\pi_\iota,\Hil_\iota,\xi_\iota)=\GNS(A_\iota,\phi_\iota)$ and
$(\pi,\Hil,\xi)=\GNS(\Afr,\phi)$.
From the free product construction we have that
$$ \Hil=\Cpx\xi\oplus
\bigoplus\Sb n\ge1 \\ \iota_1\neq\iota_2\neq\cdots\neq\iota_n \endSb
\Hilo_{\iota_1}\otimes\cdots\otimes\Hilo_{\iota_n}, $$
where $\Hilo_\iota=\Hil_\iota\ominus\Cpx\xi_\iota$.
By~\scite{\DykemaZZFPFDNT}{Lemma 3.2}, for each $\iota\in I$ there is a
projection $P_\iota:\Hil_\iota\to\Hil_{B_\iota}\eqdef\{\bh\mid b\in B_\iota\}$
such that for every $a\in A_\iota$,
$P_\iota\pi_\iota(a)\restrict_{\Hil_{B_\iota}}
=P_\iota\pi_\iota(E_\iota(a))\restrict_{\Hil_{B_\iota}}$.
Now $Y\eqdef\lspan\bigl(\{1\}\cup\Lambdao((\Bo_\iota)_{\iota\in I})\bigr)$ is a
dense subset of $B$, so
$$ \Hil_B\eqdef\overline{\{\bh\mid b\in B\}}=\overline{\{\bh\mid b\in Y\}}
=\Cpx\xi\oplus\bigoplus
\Sb n\ge1 \\ \iota_1\neq\iota_2\neq\cdots\neq\iota_n \endSb
\Hilo_{B_{\iota_1}}\otimes\cdots\otimes\Hilo_{B_{\iota_n}}, $$
where $\Hilo_{B_\iota}=\Hil_{B_\iota}\ominus\Cpx\xi_\iota$.
Let $P:\Hil\to\Hil_B$ be the projection onto $\Hil_B$.
Now $X\eqdef\lspan(\{1\}\cup\Lambdao((\Ao_\iota)_{\iota\in I}))$ is a dense
subset of $A$, and (from the free product construction) we see that
$$ P\pi(X)\restrict_{\Hil_B}\subseteq P\pi(B)\restrict_{\Hil_B}, $$
so by Lemma~\ProjGivesExp, there is a projection of norm~$1$, $E:A\to B$,
satisfying $\phi\circ E=\phi$ and given by
$P\pi(E(a))\restrict_{\Hil_B}=P\pi(a)\restrict_{\Hil_B}$.
From this  and the free product construction, we see that
$E(a)=E_\iota(a)$ whenever $a\in A_\iota\subseteq A$.
\QED

\proclaim{Definition \FPCondExp}\rm
The projection of norm~$1$, $E$, found in the above proposition is called the
{\it free product} of the $E_\iota$, and is denoted $E=\freeprodi E_\iota$.
\endproclaim

\proclaim{Corollary \ProjOntoFactor}
Let $A$ and $B$ be unital C$^*$--algebras with faithful states $\phi_A$ and
$\phi_B$, respectively.
Let
$$ (\Afr,\phi)=(A,\phi_A)*(B,\phi_B) $$
by the free product of C$^*$--algebras.
Then there is a projection of norm~$1$, $E$, from $\Afr$ onto $A$, such that
$\phi\circ E=\phi$.
\endproclaim
\demo{Proof}
Consider the projections of norm~$1$, $\id_A:A\to A$ and
$\phi_B:B\to\Cpx1\subseteq B$.
Let $E=\id_A*\phi_B$ be their free product.
\QED

\vskip3ex
\heading\S\PIFP.  Some free product C$^*$--algebras.
\endheading
\vskip3ex

\proclaim{Theorem \mostgeneral}
Let $A$ and $B$ be C$^*$--algebras with faithful states $\phi_A$, respectively,
$\phi_B$.
Consider the C$^*$--algebra reduced free product
$$ (\Afr,\phi)=(A,\phi_A)*(B,\phi_B). $$
Suppose there is a partial isometry $v\in A$, whose range projection, $q=vv^*$,
and domain projection, $p=v^*v$, are orthogonal, and such
that, for some $0<\lambda<1$, $v$ is in the spectral subspace
of $\phi_A$ associated to $\lambda^{-1}$, namely, that
$\phi_A(xv)=\lambda^{-1}\phi_A(vx)$ for every $x\in A$.
Note this implies $\phi(q)=\lambda\phi(p)<\phi(p)$.
Let
$$ A_{00}=\Cpx p+\Cpx q+(1-p-q)A(1-p-q) $$
and let $\Afr_{00}=C^*(A_{00}\cup B).$
Suppose that $q$ is equivalent in the centralizer of the restriction of $\phi$
to $\Afr_{00}$ to a subprojection of $p$,
and that the centralizer of the restriction of $\phi$ to $q\Afr_{00}q$ contains
an abelian subalgebra on which $\phi$ is diffuse
(i.e\. an abelian subalgebra to which the restriction of $\phi$ is given by
an atomless measure --- see~\scite{\DykemaZZSimplicity}{2.1}).
Suppose also that $p+q$ is full in $\Afr$.

Then $\Afr$ is simple and purely infinite.
\endproclaim
\demo{Proof}
Since $p$ is full in $\Afr$, in order to show that $\Afr$ is simple and purely
infinite, it will suffice to show that $p\Afr p$ is simple and purely
infinite.
By~\cite{\DykemaZZFaithful}, $\phi$ is faithful on $\Afr$.

By assumption there is a partial isometry, $y$, in the centralizer of
$\phi\restrict_{\Afr_{00}}$,
such that $yy^*=q$ and $p_1\eqdef y^*y\le p$.
Also by assumption, the centralizer of the restriction of $\phi$ to
$p_1\Afr_{00}p_1$ contains a diffuse abelian subalgebra.
Let $w=y^*v$.
Then $w^*w=p$ and $ww^*=p_1\le p$.
Let $p_0=p$ and for $n\ge1$ let $p_n=w^n(w^*)^n$.
Then $w$ belongs to the spectral subspace of $\Afr$ associated to
$\lambda^{-1}$, and hence
$$ \forall n\ge0\qquad\phi(p_n)=\lambda^n\phi(p). \tag{\phiofpn} $$

Let
$$ A_0=pAp+(1-p)A(1-p). $$
Then $A$ is generated by $A_0\cup\{v\}$, because
$$ A=pAp+pA(1-p)+(1-p)Ap+(1-p)A(1-p) $$
and
$$ pA(1-p)=v^*vA(1-p)=v^*qA(1-p)\subseteq v^*(1-p)A(1-p). $$
Let $\Afr_0=C^*(A_0\cup B)$.
Then $\Afr=C^*(\Afr_0\cup\{w\})$, and thus
$$ p\Afr p=C^*(p\Afr_0p\cup\{w\}). $$
In $p\Afr p$, $w$ is a proper isometry which, as we will see, is loosely
speaking as free as it can be from $p\Afr_0p$, with amalgamation over $pA_0p$.

Throughout the proof we will use a projection of norm~$1$,
$E=E^\Afr_{A_0}:\Afr\to A_0$,
which is hereby defined to be $E^\Afr_{A_0}=E^A_{A_0}\circ E^\Afr_A$, where
$E^\Afr_A:\Afr\to A$ is the conditional expectation onto $A$ obtained from
Corollary~\ProjOntoFactor, and where $E^A_{A_0}:A\to A_0$ is the
conditional expectation onto $A_0$ defined by
$$ E^A_{A_0}(a)=pap+(1-p)a(1-p). $$
Note that $\phi\circ E=\phi$.

Let $\Theta$ be the set of all
$$ x=x_1x_2\cdots x_n\in\Lambdao\bigl((p\Afr_0p)\oup,
\{w^k\mid k\ge1\}\cup\{(w^*)^k\mid k\ge1\}\bigr) \tag{\xinTheta} $$
such that whenever $2\le j\le n-1$ and $x_j\in(p\Afr_0p)\oup$,
$$ \align
\text{if }x_{j-1}=w\text{ and }x_{j+1}=w^*\qquad&\text{then }
x_j\in p\Afr_0p\ominus pA_0p \\
\text{if }x_{j-1}=w^*\text{ and }x_{j+1}=w\qquad&\text{then }
x_j\in p_1\Afr_0p_1\ominus\bigl(y^*(qA_0q)y\bigr).
\endalign $$
The restriction of $E$ to $p\Afr_0p$ provides a projection of norm~$1$ onto
$pA_0p$, and $y^*E(y\cdot y^*)y$ provides a projection of norm$~1$ from
$p_1\Afr_0p_1$ onto $y^*(qA_0q)y$, so Lemma~\Eominus{} applies.
Notice that
$$ \align
waw^*&=wpapw^* \\
wA_0w^*&\subseteq p\Afr_0p \\
w^*aw&=w^*p_1ap_1w, \\
w^*y^*(qA_0q)yw&\subseteq p\Afr_0p
\endalign $$
and hence
$$ \lspan(\{1\}\cup\Theta)=\lspan\Bigl(\{1\}\cup
\Lambdao\bigl((p\Afr_0p)\oup,
\{w^k\mid k\ge1\}\cup\{(w^*)^k\mid k\ge1\}\bigr)\Bigr) $$
is the $*$--algebra generated by $\Afr_0\cup\{w\}$.

Like in~\cite{\DykemaZZFPFDNT}, for $x\in\Theta$ let $t_j(x)$ be the number of
$w$ minus the number of $w^*$ appearing in the first $j$ letters of $x$.
Thus, if $l(x)$ is the number of letters of $x$ (i.e\. $l(x)=n$ when $x$ is as
in~(\xinTheta)), after setting $t_0(x)=0$ we thus define
$$ t_j(x)=\cases
t_{j-1}(x)&\text{ if the $j$th letter of $x$ is from }\Afr_0\oup \\
t_{j-1}(x)+k&\text{ if the $j$th letter of $x$ is }w^k \\
t_{j-1}(x)-k&\text{ if the $j$th letter of $x$ is }(w^*)^k, \\
\endcases $$
for each $1\le j\le l(x)$.
We will use interval notation to denote subsets of the integers.
Thus, for example,
$$ \gather
[0,n]\text{ will mean }\{0,1,2,\ldots,n\} \\
[0,\infty)\text{ will mean }\{0,1,2,\ldots\} \\
(-\infty,0]\text{ will mean }\{\ldots,-2,-1,0\} \\
(-\infty,\infty)\text{ will mean }\Integers.
\endgather $$
For every interval, $I$, of $\Integers$ which contains $0$, let
$$ \Theta_I=\{x\in\Theta\mid t_{l(x)}=0
\text{ and }\forall1\le j\le l(x),\,t_j(x)\in I\}. $$
Then $\lspan(\Theta_I\cup\{1\})$ is a $*$--subalgebra of $p\Afr p$.
Let $\Afr_I=\clspan(\Theta_I\cup\{1\})$.

There is an injective endomorphism, $\sigma$, of $\Afrinf$ given by
$\sigma(a)=waw^*$.
Since $p\Afr p=C^*(\Afrinf\cup\{w\})$, it follows that $p\Afr p$ is a
quotient of the (universal) C$^*$--algebra crossed product
$$ \Afrinf\rtimes_\sigma\Naturals. $$
We will use~\scite{\DykemaRordamZZPI}{Theorem 2.1(ii)} to prove that
$\Afrinf\rtimes_\sigma\Naturals$ is
simple and purely infinite, which will imply that $p\Afr p$ is simple and
purely infinite.
In particular, Claim~\alphaouter{} and Claim~\HeredProj{} below will
show that the endomorphism $\sigma$ satisfies the hypotheses
of~\scite{\DykemaRordamZZPI}{Theorem 2.1(ii)}.

Let $\Afrtinf$ be the C$^*$--algebra inductive limit
$$ \Afrinf\overset\sigma\to\to\Afrinf\overset\sigma\to\to
\Afrinf\overset\sigma\to\to\cdots\to\Afrtinf $$
and for $n\ge 1$ consider the defining $*$--homomorphisms
$\mu_n:\Afrinf\to\Afrtinf$ such that $\mu_{n+1}\circ\sigma=\mu_n$.
There is an automorphism $\alpha$ of $\Afrtinf$ defined by
$\alpha(\mu_n(a))=\mu_n(\sigma(a))$.

\proclaim{Claim \alphaouter}
For each $m\ge1$, the automorphism $\alpha^m$ of $\Afrtinf$ is outer.
\endproclaim
\demo{Proof}
If $\alpha^m$ is inner than $\alpha^m(a)=a$ for some $a\in\Afrtinf$, $a\ge0$,
$a\ne0$.
Let $p_{k-l}=\mu_l(\sigma^k(p))$.
This coincides with the old definition of $p_{k-l}$ if $k-l\ge0$, and $p_n$ is
an approximate identity for $\Afrtinf$ as $n\to-\infty$.
Thus $\nm{a-p_nap_n}$ can be made arbitrarily small by choosing $n$ large and
negative.
But $\nm{a-p_nap_n}=\nm{\alpha^{km}(a-p_nap_n)}=\nm{a-p_{n+km}ap_{n+km}}$.
Since $\phi(p_{n+km}ap_{n+km})\le\nm a\phi(p_{n+km})$ and since by~(\phiofpn)
$\phi(p_{n+km})$ tends to $0$ as $k\to\infty$, we may conclude that
$\phi(a)=0$.
But since $\phi$ is faithful, this implies $a=0$, which contradicts the choice
of $a$.

This completes the proof of Claim~\alphaouter.
\enddemo

\proclaim{Claim \phiThetao}
If $x\in\Thetainf$ then $\phi(x)=0$.
\endproclaim
\demo{Proof}
Let $x=x_1x_2\cdots x_n\in\Thetainf$.
Rewrite each $w$ appearing in $x$ as $vy^*$ and each $w^*$ as $yv^*$.
Now group together all occurances of $y^*$, $y$ and letters from $(p\Afr_0p)\oup$
that are neighbors.
The resulting object is either an element of $(p\Afr_0p)\oup$ or is equal to
$$ u=u_1u_2\cdots u_m\in\Lambdao(\Afr_0,\{v,v^*\}), $$
where whenever $2\le j\le n-1$ and $u_j\in\Afr_0$,
$$ \align
\text{if }u_{j-1}=v\text{ and }u_{j+1}=v\qquad
&\text{then }u_j\in p\Afr_0q \\
\text{if }u_{j-1}=v\text{ and }u_{j+1}=v^*\qquad
&\text{then }u_j\in p\Afr_0p\ominus pA_0p \\
\text{if }u_{j-1}=v^*\text{ and }u_{j+1}=v\qquad
&\text{then }u_j\in q\Afr_0q\ominus qA_0q \\
\text{if }u_{j-1}=v^*\text{ and }u_{j+1}=v^*\qquad
&\text{then }u_j\in q\Afr_0p.
\endalign $$
But
$$ \Afr_0\subseteq\clspan\bigl(\{1\}\cup\Lambdao(A_0\oup,B\oup)\bigr) $$
and using the conditional expectation from $\Afr_0$ onto $A_0$, we see that
$$ p\Afr_0q\cup(p\Afr_0p\ominus pA_0p)\cup(q\Afr_0q\ominus qA_0q)\cup q\Afr_0p
\subseteq\clspan\bigl(\Lambdao(A_0\oup,B\oup)\backslash A_0\oup\bigr). $$
Therefore, it will suffice to show that $\phi(z)=0$ for every
$$ z=z_1z_2\cdots z_s\in\Lambdao\bigl(\Lambdao(A_0\oup,B\oup),\{v,v^*\}\bigr)
$$
which has the property that if
$z_j\in\Lambdao(A_0\oup,B\oup)$ and $2\le j\le s-1$ then
$$ z_j\in\Lambdao(A_0\oup,B\oup)\backslash A_0\oup. $$
But since $A_0vA_0\subseteq\ker\phi_A$, we see that
$z\in\Lambdao(A\oup,B\oup)$, and it then follows
from the freeness of $A$ and $B$ that $\phi(z)=0$.

This completes the proof of Claim~\phiThetao.
\enddemo

\proclaim{Claim \AfrFreenegone}
The subalgebras $w^*\Afr_{(-\infty,0]}w$ and $\Afr_{[0,\infty)}$ are free
with amalgamation over $pA_0p$, (with respect to the restrictions of the
conditional expectation $E$).
\endproclaim
\demo{Proof}
Since $pA_0p=w^*y^*(qA_0q)yw$ we see that
$pA_0p\subseteq w^*\Afr_{(-\infty,0]}w$, and clearly
$pA_0p\subseteq p\Afr_{[0,\infty)}p$.
In order to show freeness with amalgamation, (and refering to Lemma~\Eominus),
it will suffice to show that $E(x)=0$ whenever
$$ x\in\Lambdao(\Afr_{[0,\infty)}\ominus pA_0p,
w^*\Afr_{(-\infty,0]}w\ominus pA_0p). \tag{\heresx} $$
Let $\Theta_{(-\infty,0^-]}$ be
the set of all elements of $\Theta_{(-\infty,0]}$ which are words that begin
with $w^*$ and end with $w$.
We will show that
$$ w^*\Afr_{(-\infty,0]}w\ominus pA_0p\subseteq\clspan\Theta_{(-\infty,0^-]}.
\tag{\inThetaneg} $$
Firstly, note that
$$
w^*(\Theta_{(-\infty,0]}\backslash(p\Afr_0p))w\subseteq\Theta_{(-\infty,0^-]}.
$$
Now it will be enough to show that
$$ w^*\Afr_0w\ominus pA_0p\subseteq\Theta_{(-\infty,0^-]}. $$
But
$$ \align
w^*\Afr_0w\ominus pA_0p&=w^*(p_1\Afr_0p_1\ominus wA_0w^*)w \\
&=w^*(p_1\Afr_0p_1\ominus y^*vA_0v^*y)w \\
&=w^*(p_1\Afr_0p_1\ominus y^*(qA_0q)y)w\subseteq\Theta_{(-\infty,0^-]}.
\endalign $$
We also see that
$$ \Afr_{[0,\infty)}\ominus pA_0p\subseteq
\clspan\bigl((\Theta_{[0,\infty)}\backslash(p\Afr_0p))
\cup(p\Afr_0p\ominus pA_0p)\bigr). $$
Hence, given $x$ as in~(\heresx), in order to $E(x)=0$ we may assume
without loss of generality that
$$ x\in\Lambdao\bigl(\Theta_{(-\infty,0^-]},
(\Theta_{[0,\infty)}\backslash(p\Afr_0p))\cup(p\Afr_0p\ominus pA_0p)\bigr). $$
But then clearly
$$ x\in(\Theta_{(-\infty,\infty)}\backslash\Afr_0)\cup(p\Afr_0p\ominus pA_0p).
$$
If $x\in p\Afr_0p\ominus pA_0p$ then by Lemma~\Eominus{} $E(x)=0$.
Furthermore, using Claim~\phiThetao, that
$$ p\Afr_0p(\Theta_{(-\infty,\infty)}\backslash\Afr_0)\subseteq
\lspan(\Theta_{(-\infty,\infty)}\backslash\Afr_0) $$and Lemma~\Eominus,
we see that if $x\in\Theta_{(-\infty,\infty)}\backslash\Afr_0$ then
$E(x)=0$.

This completes the proof of
Claim~\AfrFreenegone.
\enddemo

\proclaim{Claim \AfrZeroSimple}
The C$^*$--algebra $\Afr_{(-\infty,0]}$ is simple.
\endproclaim
\demo{Proof}
$\Afr_{(-\infty,0]}$ is generated by $w^*\Afr_{(-\infty,0]}w$ and $p\Afr_0p$,
which by Claim~\AfrFreenegone{} are free with amalgamation over $pA_0p$.
Now letting $\Afr_0'$ be the
C$^*$--algebra generated by $B\cup(\Cpx p+(1-p)A(1-p))$,
using~\scite{\DykemaZZSimplicity}{2.8} we see that
$p\Afr_0p$ is generated by $pA_0p$ and $p\Afr_0'p$, which are free with respect
to the state $\phi$ (after rescaling).
Hence we see that $\Afr_{(-\infty,0]}$ is generated by $w^*\Afr_{(-\infty,0]}w$
and $p\Afr_0'p$ which are free (with amalgamation over the scalars $\Cpx p$).
But $w^*\Afr_{(-\infty,0]}w\supseteq w^*(p_1\Afr_{00}p_1)w$ and the
restriction of $\phi$ to $w^*(p_1\Afr_{00}p_1)w$ is just a rescaling
of the restriction of $\phi$ to $p_1\Afr_{00}p_1$.
We saw earlier that the centralizer of the restriction of $\phi$ to
$p_1\Afr_{00}p_1$ has an abelian
subalgebra which is diffuse with respect to $\phi$, and hence so does
the centralizer of the restriction of $\phi$ to $w^*(p_1\Afr_{00}p_1)w$.
Clearly $p\Afr_0'p\neq\Cpx$, so by~\scite{\DykemaZZSimplicity}{3.2},
$\Afr_{(-\infty,0]}$ is simple.

Hence Claim~\AfrZeroSimple{} is proved.
\enddemo

\proclaim{Claim \AfrNegSimple}
For every $n\ge0$, the C$^*$--algebra $\Afr_{(-\infty,n]}$ is simple.
\endproclaim
\demo{Proof}
We use induction on $n$.
The case $n=0$ holds by the previous claim.
Assume $n\ge1$.
Now since $\Afr_{(-\infty,0]}$ is simple, $p_1$ is full in
$\Afr_{(-\infty,0]}$, hence is full in $\Afr_{(-\infty,n]}$.
Therefore it will suffice to show that $p_1\Afr_{(-\infty,n]}p_1$ is simple.
But
$$ p_1\Afr_{(-\infty,n]}p_1=ww^*\Afr_{(-\infty,n]}ww^*\cong
w^*\Afr_{(-\infty,n]}w=\Afr_{(-\infty,n-1]}, $$
which is simple by inductive hypothesis.

Hence Claim~\AfrNegSimple{} is proved.
\enddemo

\proclaim{Claim \pnkFree}
Let $n\ge0$ and $k\ge1$ be integers.
Then $p_{n+1}\Afr_{(-\infty,n]}p_{n+1}$ and $\{p_{n+k}\}$ are free (with
amalgamation over the scalars $\Cpx p_{n+1}$) with respect
to the state $\phi$ (after rescaling).
\endproclaim
\demo{Proof}
The map $x\mapsto(w^*)^{n+1}xw^{n+1}$ is an isomorphism from
$p_{n+1}\Afr_{(-\infty,n]}p_{n+1}$ onto
$w^*\Afr_{(-\infty,0]}w$ which scales
the state $\phi$ and which sends $p_{n+k}$ to $p_{k-1}$.
Hence it will suffice to show that $w^*\Afr_{(-\infty,0]}w$ and $\{p_{k-1}\}$
are free (with amalgamation over the scalars $\Cpx p$) with respect to $\phi$
(after rescaling).
In light of Claim~\AfrFreenegone, for this it will suffice to show that
$$ E(p_{k-1})=\frac{\phi(p_{k-1})}{\phi(p)}p. \tag{\Eofpm} $$
However, $p_{k-1}\in C^*(B\cup(1-p-q)A(1-p-q)\cup\{v\})$, and
$C^*((1-p-q)A(1-p-q)\cup\{v\})$ and
$B$ are free with respect to $\phi$.
Therefore,
$$ E(p_{k-1})\in E(C^*(B\cup(1-p-q)A(1-p-q)\cup\{v\}))
=\Cpx p+\Cpx q+(1-p-q)A(1-p-q). $$
Now $p_{k-1}\le p$ so $E(p_{k-1})\in\Cpx p$.
But $E$ preserves the state $\phi$, so~(\Eofpm) follows.

This completes the proof of Claim~\pnkFree.
\enddemo

For the next claim, we will make use of the 
comparison theory for positive elements in a
C$^*$--algebra that was introduced by
J\. Cuntz~\cite{\CuntzZZAddMult},~\cite{\CuntzZZDimFunct} (see
also~\cite{\RordamZZUHFII}).
Recall that for positive elements, $a$ and $b$ of $\Afr$, Cuntz defined
$a\lesssim b$ if there are $x_j\in\Afr$ such that
$\lim_{j\to\infty}x_j^*bx_j=a$.
Recall also that $\lesssim$ is a transitive relation.

\proclaim{Claim \HeredProj}
Let $D$ be a nonzero, hereditary C$^*$--subalgebra of $\Afrinf$.
Then there is a projection in $D$ that is equivalent in
$\Afrinf$ to $p_n$ for some $n$.
\endproclaim
\demo{Proof}
Let $h\in D$, $h\ge0$, $\nm h=1$.
Since $\bigcup_{n\ge1}\Afr_{(-\infty,n]}$ is dense in $\Afrinf$, for every
$\epsilon>0$ there is $n\in\Naturals$ and $h_n\in\Afr_{(-\infty,n]}$ such that
$\nm{h_n}=1$ and $\nm{h-h_n}<\epsilon$.
Take $\epsilon<1$ and
let $f:[0,1]\to[0,1]$ be monotone increasing such that $f(1-\epsilon)=0$ and
$f(1)=1$, and let $b=f(h_n)$.
Then $b\ge0$, $b\ne0$.
Since, by Claim~\AfrNegSimple, $\Afr_{(-\infty,n]}$ is simple,
$$ b\bigl(\Afr_{(-\infty,n]}\bigr)p_{n+1}\ne\{0\}. $$
Therefore, there is $a\in p_{n+1}\Afr_{(-\infty,n]}p_{n+1}$ $a\ge0$, $a\ne0$
such that $a\lesssim b$.
By Claim~\pnkFree, $a$ and $p_{n+k}$ are free for every $k\ge1$.
Since $\lim_{k\to\infty}\phi(p_{n+k})=0$,
using~\scite{\DykemaRordamZZProj}{5.3}, we
see that for sufficiently large $k$ we have $p_{n+k}\lesssim a$.
Thus $p_{n+k}\lesssim b$ and hence there is $x\in\Afrinf$ such that
$\nm{x^*bx-p_{n+k}}<\epsilon$.
But
$$ \gather
x^*bx\ge x^*b^{\frac12}h_nb^{\frac12}x\ge(1-\epsilon)x^*bx \\
\nm{x^*bx-x^*b^{\frac12}h_nb^{\frac12}x}
 \le\epsilon\nm{x^*bx}<\epsilon(1+\epsilon) \\
\nm{x^*b^{\frac12}h_nb^{\frac12}x-x^*b^{\frac12}hb^{\frac12}x}
 \le\nm{h_n-h}\;\nm{x^*bx}<\epsilon(1+\epsilon) \\
\nm{x^*b^{\frac12}hb^{\frac12}x-p_{n+k}}<\epsilon(3+2\epsilon).
\endgather $$
By standard arguments, taking $\epsilon$ small enough we find a projection in
$\overline{h\Afrinf h}$ that is equivalent to $p_{n+k}$.

Hence Claim~\HeredProj{} is proved.
\enddemo

Now using Claim~\alphaouter{}, Claim~\HeredProj{}
and~\scite{\DykemaRordamZZPI}{Theorem 2.1(ii)} completes the proof of the
proposition.
\QED

Now we give a list (by no means complete) of examples where the above
proposition can be applied.
For $0<\lambda\le1$, with the symbol $\psi_\lambda$ we denote the state on
$M_2(\Cpx)$,
$$ \psi_\lambda(\cdot)=\Tr_2\left(\cdot\left(
\matrix\frac1{1+\lambda}&0\\0&\frac\lambda{1+\lambda}\endmatrix
\right)\right). $$

\proclaim{Examples \examplesPI}
In each of the following cases, Theorem~\mostgeneral{} applies, showing
that the reduced free product C$^*$--algebra $\Afr$ is simple and purely
infinite.
\roster
\item"(i)"
Let $0<\lambda<\mu\le1$ and let
$$ (\Afr,\phi)=(M_2(\Cpx),\psi_\lambda)*(M_2(\Cpx),\psi_\mu). $$
\item"(ii)"
Let $0<\lambda<1$ and let
$$ (\Afr,\phi)=(M_2(\Cpx),\psi_\lambda)*(C([0,1]),\tint\cdot\dif t). $$
\item"(ii)"
Let
$$ (\Afr,\phi)=(A_1\otimes F,\phi_{A_1}\otimes\phi_F)*(B,\phi_B) $$
where $A_1$ is any C$^*$--algebra, $\phi_{A_1}$ is any faithful state on $A_1$,
$F$ is a finite dimensional C$^*$--algebra, $\phi_F$ is a nontracial, faithful
state on $F$ and the centralizer of $\phi_B$ contains a unital, abelian
C$^*$--subalgebra on
which $\phi_B$ is diffuse.
\endroster
\endproclaim
\demo{Proof}
We first consider case~(i).
In the first copy of $M_2(\Cpx)$, let
$p=\left(\smallmatrix1&0\\0&0\endsmallmatrix\right)$,
$q=\left(\smallmatrix0&0\\0&1\endsmallmatrix\right)$ and
$v=\left(\smallmatrix0&0\\1&0\endsmallmatrix\right)$.
Let $p'$, $q'$ and $v'$ be the same but in the second copy of $M_2(\Cpx)$.
Then $v$ is a partial isometry as required in Theorem~\mostgeneral.
Moreover, $p$, $q$, $p'$ and $q'$ are all in the centralizer of $\phi$ and the
C$^*$--algebra generated by $\{p,q,p',q'\}$ is isomorphic to the free product,
(in the notation of~\cite{\DykemaZZSimplicity}),
$$ (\smdp\Cpx{\frac1{1+\lambda}}p\oplus\smdp\Cpx{\frac\lambda{1+\lambda}}q)
*(\smdp\Cpx{\frac1{1+\mu}}{p'}\oplus\smdp\Cpx{\frac\mu{1+\mu}}{q'}), $$
which can be described by refering to~\scite{\DykemaZZSimplicity}{2.7}.
Thus
$$ C^*(\{p,q,p',q'\})
\cong\smdp\Cpx{\frac1{1+\lambda}-\frac1{1+\mu}}{p\wedge q'}
\oplus C([a,b],M_2(\Cpx))
\oplus\smdp\Cpx{\frac1{1+\lambda}-\frac\mu{1+\mu}}{p\wedge p'}, $$
for some $0<a<b<1$, where the trace on $C([a,b],M_2(\Cpx))$ is induced by a
measure on $[a,b]$ having no atoms, and with
$$ \align
p&=1\oplus\left(\smallmatrix1&0\\0&0\endsmallmatrix\right)\oplus1 \\
q&=0\oplus\left(\smallmatrix0&0\\0&1\endsmallmatrix\right)\oplus0.
\endalign $$
Hence $q$ is equivalent to a subprojection of $p$ in 
$C^*(\{p,q,p',q'\})$ and $qC^*(\{p,q,p',q'\})q$ contains a diffuse abelian
subalgebra, so the hypotheses of Theorem~\mostgeneral{} are fulfilled.

It is clear that case~(ii) follows from case~(iii).
In order to prove that Theorem~\mostgeneral{} applies in case~(iii),
note that since $\phi_F$ is nontracial, there is a partial isometry, $v\in F$,
with $p\eqdef v^*v$ and $q\eqdef vv^*$ orthogonal and minimal projections in
$F$ such that $\phi(q)<\phi(p)$.
Moreover, $v$ is in the spectral subspace of $\phi_A\otimes\phi_F$ associated to
$\lambda^{-1}=\phi(p)/\phi(q)$.
Let $D$ be a unital abelian C$^*$--subalgebra of the centralizer of $\phi_B$ on
which $\phi_B$ is diffuse.
Consider $C^*(\{p,q\}\cup D)$.
This is isomorphic to the reduced free product of abelian C$^*$--algebras
$$ (\smp\Cpx p\oplus\smp\Cpx q\oplus\smp\Cpx{1-p-q})*(D,\phi\restrict_D) $$
or to
$$ (\smp\Cpx p\oplus\smp\Cpx q)*(D,\phi\restrict_D), $$
depending on whether $p+q=1$ or not.
By~\scite{\DykemaZZSimplicity}{5.3}, $C^*(\{p,q\}\cup D)$ is simple and hence
$p+q$ is full in it.
By~\scite{\DykemaHaagerupRordam}{4.6(i)}, there is a unitary $u\in D$ such that
$\phi(u)=0$.
Then $p$ and $u^*pu$ are free.
If $p+q=1$ then $C^*(\{1,p,u^*pu\})$ is isomorphic to the free product
$$ (\smdp\Cpx{\frac1{1+\lambda}}p\oplus\smdp\Cpx{\frac\lambda{1+\lambda}}q)
*(\smdp\Cpx{\frac1{1+\lambda}}{u^*pu}
\oplus\smdp\Cpx{\frac\lambda{1+\lambda}}{u^*qu}). $$
Refering again to~\scite{\DykemaZZSimplicity}{2.7} we see that
$$ C^*(\{p,u^*pu\})
\cong\{f:[0,b]\to M_2(\Cpx)\mid f\text{ continuous and $f(0)$ diagonal }\}
\oplus\smdp\Cpx{\frac{1-\lambda}{1+\lambda}}{p\wedge u^*pu}, $$
for some $0<b<1$,where the trace corresponds to a measure on $[0,b]$ having no
atoms, and with
$$ \align
p&=\left(\smallmatrix1&0\\0&0\endsmallmatrix\right)\oplus1 \\
q&=\left(\smallmatrix0&0\\0&1\endsmallmatrix\right)\oplus0 \\
u^*pu&=\left(
\smallmatrix t&\sqrt{t(1-t)}\\\sqrt{t(1-t)}&1-t\endsmallmatrix\right)
\oplus1 \\
u^*qu&=\left(
\smallmatrix 1-t&-\sqrt{t(1-t)}\\-\sqrt{t(1-t)}&t\endsmallmatrix\right)
\oplus0.
\endalign $$
\line{Thus $qC^*(\{p,u^*pu\})q$ contains a diffuse abelian subalgebra and $u^*qu$ is
equivalent in}
$C^*(\{p,u^*pu\})$ to a subprojection of $p$, hence $q$ is
equivalent in $C^*(\{p,q\}\cup D)$ to a subprojection of $p$.
So in the case $p+q=1$ we are done.
But if $p+q\neq1$ then by~\scite{\DykemaZZSimplicity}{2.8},
$(p+q)C^*(\{p,q\}\cup D)(p+q)$ is isomorphic to the free product of
$(\Cpx p+\Cpx q)$ and $(p+q)C^*(\{p+q\}\cup D)(p+q)$, while
by~\scite{\DykemaZZSimplicity}{3.5}, the latter algebra has a diffuse abelian
subalgebra.
Therefore, the analysis we just did for the case $p+q=1$ applies to the free
product of $(\Cpx p+\Cpx q)$ and $(p+q)C^*(\{p+q\}\cup D)(p+q)$, showing that
$qC^*(\{p,u^*pu\})q$ contains a diffuse abelian subalgebra and
$q$ is equivalent in $C^*(\{p,q\}\cup D)$ to a subprojection of $p$.
\QED

\enddocument